# On Edgeworth Expansions in Generalized Urn Models

Sh. M. Mirakhmedov,

Institute of Mathematics and Information Technologies, Uzbekistan
(E- mail: shmirakhmedov@yahoo.com)

S. Rao Jammalamadaka,

University of California, Santa Barbara, USA.

( rao@pstat.ucsb.edu )

and

Ibrahim B. Mohamed,

University of Malaya, Malaysia,

( imohamed@um.edu.my )

**Abstract.** The random vector of frequencies in a generalized urn model can be viewed as conditionally independent random variables, given their sum. Such a representation is exploited here to derive Edgeworth expansions for a "sum of functions of such frequencies", which are also called "decomposable statistics." Applying these results to urn models such as with- and without- replacement sampling schemes as well as the multicolor Pólya- Egenberger model, new results are obtained for the chi-square statistic, for the sample sum in a without replacement scheme, and for the so-called Dixon statistic that is useful in comparing 2 samples.

## 1. Introduction.

Many combinatorial problems in probability and statistics can be formulated and indeed better understood by using appropriate urn models, which are also known as random allocation schemes. Such models naturally arise in statistical mechanics, clinical trials, cryptography etc. The properties of several types of urn models have been extensively studied in both probability and statistics literature; see e.g. the books by Johnson and Kotz (1977), Kolchin et al (1978), and survey papers by Ivanov et al. (1985), Kotz and Balakrishnan (1997).

One of the more common urn models is the sampling scheme with-replacement from a finite population which contains $N$ objects, labeled 1 through $N$; the probability that the m-th object will be selected in each of the sampling steps is equal to $p_m > 0$, $p_1 + ... + p_N = 1$. If $\eta_m$ stands for the frequency of the *m*-th object in a sample of size *n* (i.e. after *n* independent selections), then the random vector (r.vec.) $(\eta_1,...,\eta_N)$ has a multinomial distribution with parameters $(n, p_1,..., p_N)$. As is well known, one important and useful property of such a multinomial r.vec. is that its distribution can be represented as the joint conditional distribution of independent random variables $\xi_1,...,\xi_N$ given their sum $\xi_1+...+\xi_N = n$, where $\xi_m$ is Poisson $(\upsilon p_m)$ for an arbitrary positive real $\upsilon$. Such a conditional representation is indeed a characteristic property of many urn models, and thus the following definition includes several common urn models as special cases.

Let $\xi = (\xi_1,...,\xi_N)$ be a r.vec. with independent and non-negative integer components such that $P\{\xi_1+...+\xi_N = n\} > 0$, for a given integer $n > 1$. Also let $\eta = (\eta_1,...,\eta_N)$ be a r.vec. whose distribution is defined by

$$\mathcal{L}(\eta_1,...,\eta_N) = \mathcal{L}(\xi_1,...,\xi_N \mid \xi_1+...+\xi_N = n), \tag{1.1}$$

where $\mathcal{L}(X)$ here and in what follows, stands for the distribution of a r.vec. $X$. Note that (1.1) implies that $P\{\eta_1+...+\eta_N = n\} = 1$. The model defined in (1.1) is what we will call a "Generalized Urn Model" (GUM): when a sample of size $n$ is drawn from an urn containing $N$ types of objects and $\eta_m$ represents the number of $m$-th type of object appearing in the sample; the distribution of the r.vec. $\xi$ defines the sample scheme through (1.1). We are interested in the following general class of statistics

$$R_N(\eta) = \sum_{m=1}^{N} f_{m,N}(\eta_m), \tag{1.2}$$

where $f_{1,N}(x),...,f_{N,N}(x)$ are Borel functions defined for non-negative $x$. The functions $f_{m,N}$ can also be allowed to be random, in which case we will assume that the r.vec. $(f_{1,N}(x_N),...,f_{N,N}(x_N))$ for any collection of real non-negative $x_1,...,x_N$, does not depend on the r.vec. $\xi$. A statistic of the type (1.2) is called a "Decomposable Statistic" (DS) in the literature. For the case when the kernel functions $f_{m,N}$ are also random, the statistic (1.2) is called a "randomized DS" (see for instance Ivanov et al (1985), Mirakhmedov (1985) and Mikhaylov (1993)). Although the terminology DS is usually reserved for the special case when $f_{m,N}$ are not random, here we will use it for either of these cases. The following three special cases of the GUMs and related DSs are most common in applications.

**A. Sample scheme with replacement.** Let $\mathcal{L}(\xi_m) = Poi(\upsilon p_m)$ be a Poisson distribution with expectation $\upsilon p_m$, where $\upsilon \in (0, \infty)$ is arbitrary, $p_m > 0$, $m = 1,...,N$ and $p_1+...+p_N = 1$; then the r.vec. $\eta$ has the multinomial distribution $M(n, p_1,...,p_N)$ and we have a sample scheme with replacement. This scheme is associated with the random allocation of $n$ particles into $N$ cells: the cells are labeled 1 through $N$, particles are allocated into cells independently of each other and the probability of a particle falling into $m$ th cell is $p_m$, $m = 1,...,N$. The classical chi-square, likelihood-ratio statistic, and the empty-cells statistic are examples of the type (1.2) mentioned above.

**B. Sample scheme without replacement**. Suppose $\mathcal{L}(\xi_m) = Bi(\omega_m, \upsilon)$ is a binomial distribution with parameters $\omega_m > 0$ and arbitrary $\upsilon \in (0,1)$, $m = 1,...,N$, then the r.vec. $\eta$ has the multi-dimensional hypergeometric distribution:

$$P\{\eta_1 = k_1,...,\eta_N = k_N\} = \binom{\Omega_N}{n}^{-1} \prod_{m=1}^{N} \binom{\omega_m}{k_m},$$

where $\Omega_N = \omega_1 + ... + \omega_N$, $k_1 + ... + k_N = n$ and $0 \leq k_m \leq \omega_m$, $m = 1,...,N$. This GUM corresponds to a sampling scheme without replacement from a stratified finite population of size $\Omega_N$. For instance, the sample sum and the standard sample-based Estimate of the Population Total, are examples of DSs of the form (1.2).

**C. Multicolor Pólya-Egenberger urn model**. Let $\mathcal{L}(\xi_m) = NB(d_m, \upsilon)$ be negative binomial distribution with $d_m > 0$ and arbitrary $\upsilon \in (0,1)$, $m = 1,...,N$. Then

$$P\{\eta_1 = k_1,...,\eta_N = k_N\} = \binom{D_N + n - 1}{n}^{-1} \prod_{m=1}^{N} \binom{d_m + k_m - 1}{k_m}, \qquad (1.3)$$

where $D_N = d_1 + ... + d_N$, is the generalized Pólya-Egenberger distribution; such a specification of the GUM corresponds to the multicolor Pólya-Egenberger urn model (see e.g. Kotz and Balakrishman (1977, Chapter 40)). For example, the number of colors that appear in the sample exactly $r$ times and the number of pairs having the same color, are statistics of the type (1.2). We note that sum of functions of "spacings-frequencies" under the hypothesis of homogeneity of two samples can be formulated as a DS in this GUM; see for instance, Holst and Rao (1981), Jammalamadaka and Schweitzer (1985) for further details and important applications to testing hypotheses.

There is extensive literature on DSs, much of it related to sampling with and without replacement from a finite population. We specifically mention a few: Mirakhmedov (1996) obtains a bound for the remainder term in CLT and Cramer's type large deviation result for a special class of GUM; Mirakhmedov (1992) and Ivchenko and Mirakhmedov (1992) consider a two-term expansion with applications to some special cases of DS in a multinomial scheme under somewhat restrictive conditions; Babu and Bai (1996) obtain Edgeworth expansion for mixtures of global and local distributions -- results that can be used when the DS is a linear function of frequencies and a GUM is defined by *identically* distributed r.v.s $\xi_m$. Such results are clearly very restrictive on the parameters of the urn model and on the kernel functions $f_{m,N}$.

The aim of this paper is threefold: First, we present a general approach that allows one to obtain an Edgeworth asymptotic expansion to *any number of terms*, for the distribution of a DS in a GUM. Second, this general approach is used to extend known results for classes of DS in the 3 special cases of GUM just mentioned. Third, we illustrate these results by obtaining general Edgeworth expansions for three special and interesting cases of DS, viz.

(i) the chi-square statistic in Case A

(ii) sample-sum in a sample scheme without replacement i.e. in Case B, and

(iii) the Dixon spacings-frequencies statistic in Case C.

The chi-square statistic is considered for the case when the number of groups increases along with the sample size, a situation that has been considered by some authors, including Holst (1972), Morris (1975), Quine and Robinson (1984, 1985) and Mirakhmedov (1987). We obtain here a three-term asymptotic expansion under very general conditions on the parameters, generalizing the results in Mirakhmedov (1992) and Ivchenko and Mirakhmedov (1992). The result in (ii) improves the main results of Mirakhmedov (1983), Bloznelis (2000), as well as parts of Theorem 1 of Hu et al (2007). Asymptotic expansion for a DS in the multicolor Pólya- Egenberger urn model and for the Dixon statistic as a special case, are obtained here for the first time.

It should be remarked that although we confine our discussion in this paper to the above 3 examples of GUM and related DS for illustrative purposes as well as to keep the length of the paper reasonable, it should be mentioned that the results derived in this paper are generally applicable to any DS in other specifications of GUM, for instance to the context of a specified random forests, a random cyclic substitutions (cf. Kolchin (1985), Pavlov and Cherepanova (2002)).

The paper is organized as follows. In Section 2 we present a systematic procedure for obtaining an asymptotic expansion for the characteristic function of a DS, to terms of any order. Our general approach is based on the so-called Bartlett's type integral formula and provides a simpler and more streamlined way of obtaining higher order approximations than what previous authors have used. The main results are presented in Section 3. For the sake of completeness and to connect to Bartlett's type formula, we also present two theorems on asymptotic normality and Berry-Esseen type bounds, showing how the current formulation helps simplify similar results obtained in Mirakhmedov (1994, 1996). Applications to the special DSs (i), (ii), and (iii) are given in Section 4, while the proofs of the main results are postponed to an Appendix.

It should be mentioned that we are dealing with triangular arrays where all the parameters of a GUM vary (including the distribution of the r.v.s $\xi_m$) when both $n$ and $N$ tend to infinity, formally through a non-decreasing sequence of positive integers $\{n_v\}$, $\{N_v\}$, as $v \to \infty$; hence it is important to express the remainder terms in our asymptotic expansions which show their explicit dependence on the $n$, $N$, distributions of the r.v.s $\xi_m$ and the kernel functions $f_{m,N}$.

In what follows $c$, $C$ with or without an index, are universal positive constants which may depend on the argument and may be different at different places; all asymptotic relations and limits are considered as $n \to \infty$, and $N = N(n) \to \infty$.

## 2. Bartlett's type formula and asymptotic expansion of the characteristic function of a DS.

Define

$$A_N = \sum_{m=1}^{N} E\xi_m, \quad B_N^2 = \sum_{m=1}^{N} Var\xi_m, \quad x_N = (n - A_N)/B_N,$$

$$\Lambda_N = \sum_{m=1}^{N} Ef_m(\xi_m), \quad \gamma_N = \frac{1}{B_N^2} \sum_{m=1}^{N} cov(f_m(\xi_m), \xi_m),$$

$$g_m(y) = f_m(y) - Ef_m(\xi_m) - \gamma_N(y - E\xi_m), \quad \hat{R}_N(\eta) = \sum_{m=1}^{N} g_m(\eta_m)$$

$$\sigma_N^2 = \sum_{m=1}^{N} Var g_m(\xi_m) = \sum_{m=1}^{N} Var f_m(\xi_m) - B_N^2 \gamma_N^2. \tag{2.1}$$

Under some mild conditions, one can show that as $n \to \infty$ and $N = N(n) \to \infty$

$$ER_N(\eta) = \Lambda_N + x_N B_N \gamma_N - \frac{1 - x_N^2}{2B_N^2} \sum_{m=1}^{N} Eg_m(\xi_m)(\xi_m - E\xi_m)^2 \ 1 + o(1) \ ,$$

$$VarR_N(\eta) = \sigma_N^2 \ 1 + o(1) \ .$$

Also $\hat{R}_N(\eta) = R_N(\eta) - \Lambda_N - x_N B_N \gamma_N$ and

$$\sum_{m=1}^{N} Eg_m(\xi_m) = 0, \quad \sum_{m=1}^{N} cov \ g_m(\xi_m), \xi_m \ = 0. \tag{2.2}$$

Let $\phi$ be a measurable function such that $E|\phi(\xi_1, \xi_2, ..., \xi_N)| < \infty$ and $\zeta_N = \xi_1 + ... + \xi_N$. We have $E(\phi(\xi_1, ..., \xi_N)|\zeta_N = n) = E\phi(\eta_1, ..., \eta_N)$, because of (1.1). This together with

$$E(\phi(\xi_1, ..., \xi_N)e^{i\tau(\zeta_N - n)}) = \sum_{k=0}^{\infty} e^{i\tau(k-n)} P\{\zeta_N = k\} E(\phi(\xi_1, ..., \xi_N)|\zeta_N = k),$$

implies, by Fourier inversion

$$E\phi(\eta_1, \eta_2, ..., \eta_N) = \frac{1}{2\pi P\{\zeta_N = n\}} \int_{-\pi}^{\pi} E\phi(\xi_1, \xi_2, ..., \xi_N) \exp\{i\tau(\zeta_N - n)\} d\tau. \tag{2.3}$$

Set

$$\Psi_N(t, \tau) = \prod_{m=1}^{N} E \exp\{it\sigma_N^{-1} g_m(\xi_m) + i\tau B_N^{-1}(\xi_m - E\xi_m)\},$$

$$\Theta_N(t, x_N) = \frac{1}{\sqrt{2\pi}} \int_{-\pi B_N}^{\pi B_N} e^{-i\tau x_N} \Psi_N(t, \tau) d\tau. \tag{2.4}$$

Then Equation (2.3) together with the inversion formula for the local probability $P \ \zeta_N = n$ gives us the following Bartlett's type formula (cf. Bartlett (1938)):

$$\varphi_N(t, x_N) =: Ee^{it\sigma_N^{-1} \hat{R}_N(\eta)} = \frac{\Theta_N(t, x_N)}{\Theta_N(0, x_N)} \tag{2.5}$$

which provides the crucial formula of interest. Special formulations of this Formula (2.5) show up in literature; see e.g. Holst (1979), Quine and Robinson (1984), Mirakhmedov (1985, 1987, 1992). Also, a very special case of (2.5) is the most commonly used formula of Erdös and Renyi (1959) for investigating the sample sum in a without-replacement scheme (see e.g. Babu and Singh (1985), Zhao et al (2004) and Hu et al (2007)). Formula (2.3) is also useful in studying large deviation problems (see e.g. Mirakhmedov (1996)).

A formal construction of the asymptotic expansion for $\varphi_N(t, x_N)$ defined in (2.5), proceeds as follows: The integrand $\Psi_N(t, \tau)$ is the characteristic function (ch.f.) of the sum of $N$ independent two-dimensional r.vec.s $g_m, \xi_m$. Because of (2.2), this sum has zero expectation, a unit covariance matrix and uncorrelated components. From Bhattacharya and Rao (1976, Chapter 2) (this reference will henceforth be referred to as BR), it is well known that under suitable conditions, this ch. f. $\Psi_N(t, \tau)$ can be approximated by a power-series in $N^{-1/2}$ whose coefficients are polynomials in $t$ and $\tau$ containing the common factor $\exp\{-(t^2 + \tau^2)/2\}$. Hence the series can be integrated wrt $\tau$ over the interval $(-\infty, \infty)$. As a result of this integration, we get a power series, say $H_N(t, x_N)$, in $N^{-1/2}$. Next, we replace $\Theta_N(0, x_N)$ by its series approximation, which is $H_N(0, x_N)$. Finally, we get the asymptotic expansion of $\varphi_N(t, x_N)$ by dividing $H_N(t, x_N)$ by $H_N(0, x_N)$.

The above algorithm, although manageable, needs long and complex calculations as we show below. Assume that $E|g_m(\xi_m)|^s < \infty$ and $E|\xi_m|^s < \infty$ for some $s \geq 3$. Let $P_{m,N}(t, \tau), m = 1, 2, \ldots$, be the well-known polynomials in $t$ and $\tau$ from the theory of the asymptotic expansion of the ch.f. of the sum of independent random vectors (see (7.3), (7.6) of BR, p.52), in our case for the quantity $g_1, \xi_2 + \ldots + g_N, \xi_N$; the degree of $P_{m,N}(t, \tau)$ is $3m$ and the minimal degree is $m+2$; the coefficients of $P_{m,N}(t, \tau)$ only involve the cumulants of the r.v.s $g_1, \xi_2, \ldots, g_N, \xi_N$ of order $m+2$ and less. Define polynomials (in $t$) of $G_{k,N}(t, x_N)$ as

$$G_{k,N}(t, x_N) = \frac{e^{x_N^2/2}}{\sqrt{2\pi}} \int_{-\infty}^{\infty} P_{k,N}(t, \tau) \exp\left\{-i\tau x_N - \frac{\tau^2}{2}\right\} d\tau, \quad k = 0, 1, 2, \ldots \quad (2.6)$$

Now define $Q_{j,N}(x_N)$ from the equation

$$\sum_{k=0}^{\infty} (-1)^k \left(\sum_{v=0}^{s-3} N^{-v/2} G_{v,N}(0, x_N)\right)^k = \sum_{j=0}^{\infty} N^{-j/2} Q_{j,N}(x_N).$$

Then

$$Q_{j,N}(x_N) = j! \sum \prod_{i=1}^{s-3} \frac{1}{j_i!} G_{i,N}^{j_i}(0, x_N), \quad (2.7)$$

where the summation is over all $(s-3)$-tuples $(j_1, j_2, \ldots, j_{s-3})$ with non-negative integers $j_i$ such that $j_1 + 2j_2 + \ldots + (s-3)j_{s-3} = j$. Let

$$W_N^{(s)}(t, x_N) = \sum_{m=0}^{s-3} N^{-m/2} \sum_{v=0}^{m} G_{v,N}(t, x_N) Q_{m-v,N}(x_N). \quad (2.8)$$

Note that $G_{0,N}(t,x_N)=1, Q_{0,N}(x_N)=1$ so that $W_N^{(3)}(t,x_N)=1$. For example

$$W_N^{(5)}(t,x_N) = 1 + \frac{1}{\sqrt{N}}(G_{1,N}(t,x_N) - G_{1,N}(0,x_N))$$

$$+ \frac{1}{N}(G_{2,N}(t,x_N) - G_{2,N}(0,x_N) - G_{1,N}(0,x_N)(G_{1,N}(t,x_N) - G_{1,N}(0,x_N))). \quad (2.9)$$

In what follows, we shall use the following additional notation

$$\hat{\sigma}_N^2 = N^{-1}\sigma_N^2, \quad \hat{B}_N^2 = N^{-1}B_N^2, \quad \hat{g}_m = g_m(\xi_m)/\hat{\sigma}_N, \quad \hat{\xi}_m = (\xi_m - E\xi_m)/\hat{B}_N,$$

$$\beta_{j,N} = N^{-j/2}\sum_{m=1}^{N} E|\hat{g}_m|^j, \quad \kappa_{j,N} = N^{-j/2}\sum_{m=1}^{N} E|\hat{\xi}_m|^j,$$

$$M_N(T) = \inf_{T\leq|\tau|\leq\pi}\sum_{m=1}^{N}(1-|E\exp\{i\tau\xi_m\}|^2) \text{ if } T\leq\pi, \text{ else } M_N(T)=\infty, \quad (2.10)$$

$$\Upsilon_{s,N} = \beta_{s,N} + \kappa_{s,N} + B_N^2 \exp\left\{-\frac{1}{8}M_N\left(0.3(B_N\kappa_{3,N})^{-1}\right)\right\}, \quad T_N = \min(\beta_{3,N}^{-1}, \mathcal{E}_N^{-1}(1)),$$

$$\mathcal{E}_N(\delta) = \frac{1}{\sqrt{M_N(0.3(B_N\kappa_{2+\delta,N}^{1/\delta})^{-1})}} + \frac{\min(B_N,\sqrt{N})}{M_N(0.3(B_N\kappa_{2+\delta,N}^{1/\delta})^{-1})}, \quad 0<\delta\leq 1.$$

Throughout the paper we assume that $|x_N|\leq c$, although the method used here allows letting $x_N$ to increase at a rate of $O((\log N)^{1/2})$ (see e.g. Mirakhmedov (1994)). In the above listed three examples of GUM the parameter $\upsilon$ can be chosen such that $x_N = 0$ (see also the beginning of Sec. 4).

**Proposition 2.1.** Let $E|g_m(\xi_m)|^s < \infty$, for some $s\geq 3$, $m=1,...,N$ and $\Upsilon_{s,N}\leq 0.01$. There exist constants $c$ and $C$ such that if $|t|\leq cT_N$, then for $j=0,1$,

$$\left|\frac{\partial^j}{\partial t^j}\left(\varphi_N(t,x_N) - e^{-\frac{t^2}{2}}W_N^{(s)}(t,x_N)\right)\right| \leq Ce^{-\frac{t^2}{8}}\Upsilon_{s,N}.$$

**3. Main results.**

We use the notations of Sec.2. The following Theorems 3.1 and 3.2 follow from Theorem 1 and 2 of Mirakhmedov (1994) and are presented here for the sake of completeness, and to connect to Bartlett's type formula (2.5); also their application to DS in our examples of GUM gives weaker conditions for asymptotic normality and improved Berry-Esseen type bound than are known before.

Let $I\{A\}$ stand for the indicator function of the set $A$ and

$$\mathcal{L}_{1,N}(\varepsilon) = \frac{1}{N^{3/2}}\sum_{m=1}^{N}E|\hat{\xi}_m|^3 I\{|\hat{\xi}_m|\leq\varepsilon\}, \quad \mathcal{L}_{2,N}(\varepsilon) = \frac{1}{N}\sum_{m=1}^{N}E\hat{\xi}_m^2 I\{|\hat{\xi}_m|>\varepsilon\},$$

$$L_{2,N}(\varepsilon) = \sum_{m=1}^{N}E\hat{g}_m^2 I\{|\hat{g}_m|>\varepsilon\}.$$

**Theorem 3.1.** If for arbitrary $\varepsilon > 0$

(i) $\quad L_{2,N}(\varepsilon)\to 0$, $\quad\quad\quad\quad\quad\quad\quad\quad\quad\quad\quad\quad\quad\quad\quad\quad\quad\quad\quad\quad\quad\quad\quad\quad (3.1)$

(ii) $\mathcal{L}_{2,N}(\varepsilon) \to 0$,

(iii) $M_N(\pi(4B_N \mathcal{L}_{1,N}(\varepsilon))^{-1}) \to \infty$,

(iv) $\min(B_N, \sqrt{N}) = o(M_N(\pi(4B_N \mathcal{L}_{1,N}(\varepsilon))^{-1}))$,

then the statistic $R_N(\eta)$ has an asymptotic normal distribution with expectation $\Lambda_N + x_N B_N \gamma_N$ and variance $\sigma_N^2$, given in (2.1).

**Remark 3.1.** For all the 3 examples of GUM we consider, conditions (ii), (iii) and (iv), being conditions on the parameters of the urn model, are automatically satisfied under very general set-up (see Sec. 4), so that all we need is to check the Lindeberg's condition (i) for ensuring the asymptotic normality of the DS.

Let $E|g_m(\xi_m)|^s < \infty$, for some $s \geq 3$. Define $\mathbb{W}_N^{(s)}(u, x_N)$ so that

$$\int_{-\infty}^{\infty} e^{itu} d\mathbb{W}_N^{(s)}(u, x_N) = W_N^{(s)}(t, x_N) e^{-\frac{t^2}{2}}. \tag{3.2}$$

The function $\mathbb{W}_N^{(s)}(u, x_N)$ can be obtained by formally substituting

$$(-1)^\nu \frac{d^\nu}{du^\nu} \Phi(u) = -e^{-u^2/2} H_{\nu-1}(u)/\sqrt{2\pi}, \text{ where } \Phi(u) = \frac{1}{\sqrt{2\pi}} \int_{-\infty}^{u} e^{-\frac{t^2}{2}} dt,$$

for $(it)^\nu$ for each $\nu$ in the expression for $W_N^{(s)}(t, x_N)$ (see Lemma 7.2 of BR, p. 53), where $H_\nu(x)$ is the $\nu$-th order Hermite-Chebishev polynomial. Note that $\mathbb{W}_N^{(3)}(u, x_N) = \Phi(u)$. Set

$$\Delta_N^{(s)} = \sup_{-\infty < u < \infty} \left| P\{R_N(\eta) < u\sigma_N + \Lambda_N + x_N B_N \gamma_N\} - \mathbb{W}_N^{(s)}(u, x_N) \right|,$$

$$\chi_N(a, b) = \mathrm{I}\{a < b\} \int_{a \leq |t| \leq b} \left|\frac{\varphi_N(t, x_N)}{t}\right| dt.$$

**Theorem 3.2.** Let $0 < \delta \leq 1$. Then, there exists a constant $C$ such that

$$\Delta_N^{(3)} \leq C(\beta_{2+\delta, N} + \kappa_{2+\delta, N} + \mathcal{E}_N(\delta)).$$

**Theorem 3.3.** Let $E|g_m(\xi_m)|^s < \infty$, for some $s \geq 3$, $m = 1, 2, ..., N$. There exist constants $c$ and $C$ such that $\Delta_N^{(s)} \leq C\Upsilon_{s,N} + \chi_N(cT_N, \beta_{s,N}^{-1})$.

**Theorem 3.4.** Let the statistic $R_N(\eta)$ be a lattice r.v. with span $h$ and a set of possible values in $\Re$. If $E|g(\xi_m)|^s < \infty$, for some $s \geq 3$, $m = 1, 2, ..., N$ then there exist constants $c$ and $C$ such that uniformly in $z \in \Re$

$$\sup_{z \in \Re} \left| \frac{\sigma_N}{h} P\{R_N(\eta) = z\} - \frac{d}{du_z} \mathbb{W}_N^{(s)}(u_z, x_N) \right| \leq C\Upsilon_{s,N} + \bar{\chi}_N(cT_N, \pi\sigma_N/h),$$

where $u_z = (z - \Lambda_N - x_N B_N \gamma_N)/\sigma_N$ and

$$\bar{\chi}_N(a,b) = \mathrm{I}\{a<b\} \int_{a \leq |t| \leq b} |\varphi_N(t, x_N)| dt .$$

The following general bounds for $\chi_N(a,b)$ are useful in applications. Write

$$\psi_m(t,\tau) = E\exp\{itf_{m,N}(\xi_m) + i\tau\xi_m\}, \quad d_N(a,b) = 1 - \sup_{\substack{a\sigma_N^{-1} \leq |t| \leq b\sigma_N^{-1} \\ |\tau| \leq \pi}} N^{-1}\sum_{m=1}^{N}|\psi_m(t,\tau)|^2, \qquad (3.3)$$

$$H_m(t,\tau) = E\left\langle tf_{m,N}^*(\xi_m) + \tau\xi_m^*\right\rangle^2, \qquad \bar{H}_N(a,b) = \inf_{\substack{a\sigma_N^{-1} \leq |t| \leq b\sigma_N^{-1} \\ |\tau| \leq \pi}} \frac{1}{N}\sum_{m=1}^{N} H_m\left(\frac{t}{2\pi}, \frac{\tau}{2\pi}\right). \qquad (3.4)$$

where $\langle a \rangle$ stands for the distance between real $a$ and integers. Here and in what follows, for a given r.v. $\zeta$ we define $\zeta^* = \zeta - \zeta'$, where $\zeta'$ is an independent copy of $\zeta$. Then

$$\chi_N(a,b) \leq CB_N \ln(b\sigma_N^{-1})\exp\left\{-\frac{1}{2}Nd_N(a,b)\right\}, \qquad (3.5)$$

$$\chi_N(a,b) \leq CB_N \ln(b\sigma_N^{-1})\exp\left\{-\frac{1}{2}N\bar{H}_N(a,b)\right\} \qquad (3.6)$$

and

$$\bar{\chi}_N(a,b) \leq C\sigma_N B_N \exp\left\{-\frac{1}{2}N\bar{H}_N(a,b)\right\}. \qquad (3.7)$$

These inequalities (3.5)-(3.7) follow from the following arguments: From formula (2.3) it follows that for the ch. f. $\varphi_N(t, x_N)$ in Equation (2.5), one can write the product $\prod_{m=1}^{N}\psi_m(t,\tau)$ instead of $\Psi_N(t,\tau)$. Since $\Theta_N\left(0, x_N\right) \geq c$ (cf. (5.7) below), inequality (3.5) follows by using the fact that $x < e^{(x^2-1)/2}$. On the other hand by Lemma 4 of Mukhin (1991), we have

$$4H_m\left(\frac{t}{2\pi}, \frac{\tau}{2\pi}\right) \leq 1 - |\psi_m(t,\tau)|^2 \leq 2\pi^2 H_m\left(\frac{t}{2\pi}, \frac{\tau}{2\pi}\right).$$

This inequality together with (3.5) implies the inequalities (3.6) and (3.7).

**Remark 3.2.** A DS of the special form

$$X_N^2 = \sum_{m=1}^{N}\eta_m^2$$

arises in many problems in statistics and in discrete probability (see e.g. Sec. 4, and Pavlov and Cherepanova (2002)). This DS is a lattice r.v. with span equal to two. Also

$$H_m(t,\tau) = \sum_{k,l}\langle v_{r,l}\rangle^2 P(\xi_m = k)P(\xi_m' = l), \qquad (3.8)$$

where $v_{r,l} = (k-l)((k+l)t + \tau)$. As in Lemma 2 of Pavlov and Cherepanova (2002), one can prove that for all real $t$ and $\tau$ such that $|t| \leq 1/4, |\tau| \leq 1/2$ and any non-negative integer $k$ and $l$

$$\max\{\langle v_{k,l}\rangle, \langle v_{k+1,l}\rangle, \langle v_{k+2,l}\rangle\} \geq \frac{|t|}{2}.$$

From this it follows that if

$$\sum_{l=0}^{\infty}\sum_{j=0}^{\infty} P\{\xi_m = k_{j,l}\} P\{\xi'_m = l\} \geq c > 0, \qquad (3.9)$$

then $\bar{H}_N(a, \pi\sigma_N/2) \geq a^2/4\sigma_N^2$, where for each $l = 0, 1, 2, ...$ $k_{j,l}$ is defined such that $\max\{\langle v_{3j,l}\rangle, \langle v_{3j+1,l}\rangle, \langle v_{3j+2,l}\rangle\} = \langle v_{k_{j,l},l}\rangle$, $j = 0, 1, 2, ...$.

## 4. Applications

In what follows, we will use the notations of the preceding sections, keeping in the mind that the distribution of the r.v. $\xi_m$ is what is relevant for the particular GUM under consideration. Note that in all our examples of the GUM, the distributions of the r.v.'s $\xi_m$ depend on an arbitrary parameter $\upsilon$, which can be chosen in a suitably convenient manner. We will thus choose the parameter $\upsilon$ such that $A_N = n$, in which case $x_N = 0$, and hence the terms of asymptotic expansion, i.e. the function $\mathbb{W}_N^{(s)}(u,0)$ is considerably simplified. For example: it is known that $P_{0,N}(t,\tau) = 1$,

$$P_{1,N}(t,\tau) = \frac{i^3}{6N}\sum_{m=1}^{N} E(t\hat{g}_m + \tau\hat{\xi}_m)^3,$$

$$P_{2,N}(t,\tau) = \frac{i^4}{24N}\sum_{m=1}^{N}(E(t\hat{g}_m + \tau\hat{\xi}_m)^4 - 3(E(t\hat{g}_m + \tau\hat{\xi}_m)^2)^2) + \frac{1}{2}P_{1,N}^2(t,\tau).$$

Therefore, from (2.6), (2.9), and (3.2)

$$\mathbb{W}_N^{(5)}(u,0) = \Phi(u) - \frac{e^{-u^2/2}}{\sqrt{2\pi N}}\left(\frac{u^2-1}{6}\alpha_{3,0,N} - \frac{1}{2}\alpha_{1,2,N}\right) - \frac{e^{-u^2/2}}{\sqrt{2\pi N}}\left\{\frac{u^5 - 10u^3 + 15u}{72}\alpha_{3,0,N}^2\right.$$

$$+ \frac{u^3 - 3u}{24}\left(\alpha_{4,0,N} - \frac{3}{N}\sum_{m=1}^{N}\hat{\alpha}_{20m}^2 - 3\alpha_{2,1,N}^2 - 2\alpha_{3,0,N}\alpha_{1,2,N}\right)$$

$$+ \frac{u}{8}\left(3\alpha_{1,2,N}^2 + 2\alpha_{2,1,N}\alpha_{0,3,N} - 2\alpha_{2,2,N} + \frac{4}{N}\sum_{m=1}^{N}\hat{\alpha}_{11m}^2 + \frac{2}{N}\sum_{m=1}^{N}\hat{\alpha}_{20m}\hat{\alpha}_{02m}\right)\right\}. \qquad (4.1)$$

where $\hat{\alpha}_{ijm} = E\hat{g}_m^i\hat{\xi}_m^j$, $\alpha_{i,j,N} = N^{-1}\sum_{m=1}^{N}\hat{\alpha}_{ijm}$.

In what follows, just in order to keep our calculations simple, we will restrict ourselves to such a three-term asymptotic expansion given above.

**4.1 Example A.** The r.vec. $\eta = (\eta_1, ..., \eta_N)$ has the multinomial distribution $M(n, p_1, ..., p_N)$, $p_m > 0$, $m = 1, ..., N$, $p_1 + ... + p_N = 1$, and we take $\mathcal{L}(\xi_m) = Poi(np_m)$. We assume that $N = N(n) \to \infty$, $\max_{1 \leq m \leq N} p_m \to 0$ as $n \to \infty$. We take $\lambda = n/N$, $\lambda_m = np_m$ and $\mathrm{P}_{iN} = p_1^i + ... + p_N^i$.

In this classical scheme, since the best conditions for asymptotic normality and the Berry-Esseen type bound of DS are already given in Mirakhmedov (1996 and 2007), we concentrate our attention on the asymptotic expansion results.

**Theorem 4.1**. Let the statistic $R_N(\eta)$ be a lattice r.v. with span $h$ and a set of possible values $\Re$. If $E|g(\xi_m)|^5 < \infty$, $m = 1, 2, ..., N$ and $P_{2N} \leq (10\ln n)^{-1}$ then uniformly in $z \in \Re$

$$\left| \frac{\sigma_N}{h} P\{R_N(\eta) = z\} - \frac{d}{du_z} W_N^{(5)}(u_z, 0) \right|$$

$$\leq C \left( \beta_{5,N} + (P_{3N} + n^{-2})^{3/4} + \sigma_N \sqrt{n} \exp\left\{ -\frac{1}{2} N\bar{H}_N\left( cT_N, \frac{\pi\sigma_N}{h} \right) \right\} \right),$$

where $u_z = (z - \Lambda_N)/\sigma_N$ and $T_N$ is defined as in Sec 3 with

$$\mathcal{E}_N(1) = \sqrt{n^{-1} + P_{2N}} (1 + \min(1, \lambda^{-1/2})\sqrt{1 + nP_{2N}}).$$

The particular DS $X_N^k := \sum_{m=1}^{N} \eta_m^k$ for any integer $k > 1$, is a special case of Theorem 4.1. We shall focus on the most important application, the chi-square type statistic $X_N^2$. As stated before, $X_N^2$ is lattice with span equal two; also in this case $g_m(\xi_m) = (\xi_m^2 - \lambda_m(\lambda_m + 1)) - (2nP_{2N} + 1)(\xi_m - \lambda_m)$. Hence

$$\Lambda_N = n(1 + nP_{2N}),$$

$$\sigma_N^2 = 2n^2 P_{2N} + 4n^3(P_{3N} - P_{2N}^2) = N(2n\lambda P_{2N} + 4n^2\lambda(P_{3N} - P_{2N}^2)) := N\hat{\sigma}_N^2,$$

$$\alpha_{12N} = 2\lambda\hat{\sigma}_N^{-1} P_{2N}, \quad \alpha_{21N} = 4\sqrt{n}\lambda\hat{\sigma}_N^{-2}(P_{2N} + 12n(P_{3N} - P_{2N}^2)),$$

$$\alpha_{30N} = n\lambda\hat{\sigma}_N^{-3} \left[ 4P_{2N} + 2n(16P_{3N} - 9P_{2N}^2) + 8n^2(4P_{4N} - 9P_{2N}P_{3N} + 5P_{2N}^3) \right],$$

$$\alpha_{40N} = \hat{\sigma}_N^{-4} n\lambda \left[ 8P_{2N} + n(164P_{2N}^2 - 17P_{3N}) + n^2(636P_{4N} - 768P_{2N}P_{3N} + 192P_{2N}^3) \right.$$
$$+ n^3(448P_{5N} - 1120P_{2N}P_{4N} + 912P_{2N}^2 P_{3N} - 240P_{2N}^4)$$
$$\left. + 48n^4(P_{6N} - 4P_{2N}P_{5N} + 6P_{2N}^2 P_{4N} - 4P_{2N}^3 P_{3N} + P_{2N}^5) \right],$$

$$\alpha_{22N} = (\lambda\hat{\sigma}_N^2)^{-1}\lambda \left[ 8nP_{2N} + 2n^2(19P_{3N} - 14P_{2N}^2) + 12n^3(P_{4N} - 2P_{2N}P_{3N} + P_{2N}^3) - 1 \right],$$

$$\frac{1}{N}\sum_{m=1}^{N} \tilde{\alpha}_{20m}^2 = \hat{\sigma}_N^{-4} n\lambda \left[ 3P_{2N} - 2NP_{3N} + 4N^2 P_{4N} + 16N^3(P_{5N} - 2P_{4N}P_{2N} + P_{3N}P_{2N}) \right.$$
$$\left. + 16N^4(P_{6N} - 4P_{2N}P_{5N} + 6P_{4N}P_{2N}^2 - 4P_{3N}P_{2N}^3 + P_{2N}^5) \right],$$

$$\frac{4}{N}\sum_{m=1}^{N} \tilde{\alpha}_{11m}^2 + \frac{2}{N}\sum_{m=1}^{N} \tilde{\alpha}_{20m}\tilde{\alpha}_{02m} = (\lambda\hat{\sigma}_N^2)^{-1} \left[ 2n^2\lambda P_{3N} + 2n^3\lambda(2P_{4N} - 6P_{3N}P_{2N} + 3P_{2N}^3) \right].$$

**Corollary 4.1**. Let $c_1 \leq Np_m \leq c_2$ for some positive $c_1$, $c_2$ and all $m = 1, ..., N$ then uniformly in $b \in \{n + 2k, k = 0, 1, ..., n(n-1)/2\}$, the set of possible values of the r. v. $X_N^2$, one has

$$\left| \frac{\sigma_N}{2} P\{X_N^2 = b\} - \frac{d}{du_b} \mathbb{W}_N^{(5)}(u_b, 0) \right| \le C\left( \frac{1}{N^{3/2}} + \frac{1}{n\lambda^{3/2}} + n\lambda \exp\left\{ -\frac{cN}{\lambda \max(1,\lambda)} \right\} \right),$$

where $u_b = (b - \Lambda_N)/\sigma_N$, the exact formulae for $\Lambda_N, \sigma_N^2$ and the terms of $\mathbb{W}_N^{(5)}(u_b, 0)$ are given above.

The following Corollary 4.2 follows from Corollary 4.1 by using the Euler-Maclaurin sum formula. We state just a two-term asymptotic expansion to keep the expressions simple.

**Corollary 4.2**. Let $c_1 \le Np_m \le c_2$ for some positive $c_1, c_2$ and all $m = 1,...,N$. Then

$$\left| P\{X_N^2 < u\sigma_N + \Lambda_N\} - \Phi(u) - \frac{e^{-u^2/2}}{\sqrt{2\pi N}} \left[ \frac{1-u^2}{6} \alpha_{30N} + \frac{\lambda P_{2N}}{\hat{\sigma}_N} \right. \right.$$

$$\left. \left. + \frac{2}{\hat{\sigma}_N} S_1\left( \frac{1}{2}(u\sigma_N + \Lambda_N) \right) \right] \right| \le C\left( \frac{1}{N} + \frac{1}{n\lambda} + n\lambda \exp\left\{ -\frac{cN}{\lambda \max(1,\lambda)} \right\} \right), \qquad (4.2)$$

where $S_1(x) = x - [x] + 1/2$ is well-known periodic function of period one (see for instance BR, p.254), and comes up here due to the Euler-Maclaurin sum formula.

We may remark here that Corollary 4.2 already considerably improves Theorem 5 of Ivchenko and Mirakhmedov (1992), which states the inequality (4.2) with $\exp\{-N\lambda^l e^{-2\lambda}\}, l > 0$, instead of the exponential term, which makes sense under the additional restriction that $\lambda = O(\ln N)$.

**Remark 4.1**. Application of Theorems 3.3 and 3.4 to the log-likelihood statistic $L_N = \sum_{m=1}^{N} \eta_m \ln \eta_m$, and to the count-statistics, $\mu_r = \sum_{m=1}^{N} \mathbb{I}\{\eta_m = r\}$ give results similar to Theorems 4 and 6 respectively of Ivchenko and Mirakhmedov (1992), but one can obtain additional terms in the expansions they provide.

A DS with kernel functions $f_{m,N} = f_N$ for all $m = 1, 2, ..., N$ is called a "symmetric DS"; For example, the $X_N^2$, $\mu_r$ and $L_N$ are all symmetric DS. It is well-known (see e.g. Holst (1972), Quine and Robinson (1985), Mirakhmedov (1987)) that the chi-square test is asymptotically most powerful (AMP) within the class of symmetric tests, i.e. among tests based on symmetric DS, for testing the hypothesis of uniformity against the sequence of alternatives $H_{1n}$ given by:

$$p_m = \frac{1}{N}\left( 1 + \frac{\vartheta_m}{(n\lambda)^{1/4}} \right), \; m=1,2,...,N\,; \; \sum_{m=1}^{N} \vartheta_m = 0 \text{ and } 0 < C_1 \le \frac{1}{N}\sum_{m=1}^{N} \vartheta_m^2 \le C_2 < \infty.$$

Moreover the chi-square test is the unique AMP test for $\lambda$ bounded away from zero and infinity; on the other hand if $\lambda \to 0$ or $\lambda \to \infty$ then there exist other AMP symmetric tests, for example, the empty cells test when $\lambda \to 0$, and the log-likelihood test when $\lambda \to \infty$. In view of this, Ivchenko and Mirakhmedov (1992) introduced and studied the "second order asymptotic efficiency" (SOAE) of symmetric tests wrt the chi-square test. Investigation of the SOAE is based on the asymptotic expansion of the power function of such tests. In the case $\lambda \to 0$ they have shown that SOAE may arise only if $n = O(N^{3/4})$; for example the empty cells test based on the statistic $\mu_0$ is SOAE for this situation; for the case $\lambda \to \infty$ they could only note that when $\lambda = O(\ln N)$ the SOAE test does not exist, because of the restrictive choice of $\lambda = O(\ln N)$ needed in their asymptotic expansions. Therefore, they pointed out that the SOAE problem is open for $\lambda \to \infty$. Corollary 4.2 does resolve this problem showing that the chi-square test is still optimal in the sense of SOAE if $n = o(N^{3/2})$; it is also SOAE wrt the log-likelihood test for $n \geq N^{3/2}$ (cf. Ivchenko and Mirakhmedov(1992) for further discussion).

**4.2. Example B.** Now we consider the sample scheme without replacement from a stratified population of size $\Omega_N$; the strata are indexed by $m = 1, ..., N$; $\omega_m$ is the size of the $m$-th stratum, with $\Omega_N = \omega_1 + ... + \omega_N$; and $\eta_m$ is the number of elements of the $m$-th stratum appearing in a sample of size $n$. In this scheme $\mathcal{L}(\xi_m) = Bi(\omega_m, \upsilon)$, where $\upsilon \in (0,1)$ is arbitrary. We choose $\upsilon = p =: n/\Omega_N$, $q = 1 - p$, so that $x_N = 0$. Set $\bar{\omega}_N = \max_{1 \leq m \leq N} \omega_m$, $\Omega_{2,N} = \omega_1^2 + ... + \omega_N^2$. We consider the case where the strata sizes $\omega_m$ may increase together with $N$, but satisfy the following condition

$$\bar{\omega}_N = o((nq)^{1/4}). \tag{4.3}$$

**Theorem 4.2**. If the Lindberg's condition (3.1) is satisfied along with Condition (4.3), then as $nq \to \infty$, $R_N(\eta)$ has the asymptotic normal distribution with expectation $\Lambda_N$ and variance $\sigma_N^2$ as given in (2.1).

**Theorem 4.3**. For arbitrary $\delta \in (0,1]$ there exists a constant $C$ such that

$$\Delta_N^{(3)} \leq C \left( \beta_{2+\delta, N} + \left( \frac{\bar{\omega}_N}{nq} \right)^{\delta/2} + \frac{\bar{\omega}_N^2}{\sqrt{nq}} \right).$$

**Remark 4.2.** The term $\bar{\omega}_N^2 / \sqrt{nq}$ can be replaced by $\bar{\omega}_N \max \left| 1 - 6pq + 3nq\Omega_{2,N}\Omega_N^{-2} \right| / \sqrt{nq}$. If $\bar{\omega}_N \leq (nq)^{(1-\delta)/(4-\delta)}$ then the second term on the rhs dominates the third one.

**Theorem 4.4**. Let $E|g_m(\xi_m)|^5 < \infty$ then there exist constants $c$ and $C$ such that

$$\Delta_N^{(5)} \leq C \left( \beta_{5,N} + \left( \frac{\bar{\omega}_N}{nq} \right)^{\frac{3}{2}} + \chi_N(cT_N, \beta_{5,N}^{-1}) \right).$$

Let now the elements of $m$-th stratum be independent r.v.s $X_{m,1}, ..., X_{m,\omega_m}$, $m = 1, ..., N$. We draw a sample of size $n$ without replacement from the entire population. Define the indicator r.v.s $\eta_{mi}$ which equal one if an element $X_{mi}$ of the $m$ th stratum appears in the sample, or else it is zero, so that $\eta_{m1} + ... + \eta_{m\omega_m} = \eta_m$. Then $S_{n,N}^{(m)} = \sum_{i=1}^{\omega_m} X_{mi}\eta_{mi}$ represents the sum of elements of the $m$ th stratum which appear in the sample, and the sum of all the elements in the sample, the "sample-sum" given by $S_{n,N} = \sum_{m=1}^{N} S_{n,N}^{(m)}$, is a DS.

Assume that the r.v.s $X_{m,1},\ldots,X_{m,\omega_m}$ have a common distribution, same as that of a r.v. $Y_m$, $m = 1,\ldots,N$. We also assume that $Y_1,\ldots,Y_N$ are independent r.v.s. Then the r.v. $S_{n,N}$ is distributionally equal to a DS with $f_{mN}(0) = 0$, $f_{mN}(j) = X_{m,1} + \ldots + X_{m,j}$, $m = 1,\ldots,N$:

$$\mathcal{L}(S_{n,N}) = \mathcal{L}\left(\sum_{m=1}^{N}\left(\sum_{j=1}^{\eta_m} X_{m,j} I\{\eta_m \geq 1\}\right)\right). \tag{4.4}$$

Suppose $E|Y_m|^s < \infty$ for some $s \geq 3$. Then the expressions in (2.1) have the following form:

$$f_{mN}(\xi_m) = \sum_{j=1}^{\xi_m} X_{m,j} I\{\xi_m \geq 1\}, \tag{4.5}$$

$$g_m(\xi_m) = \sum_{j=1}^{\xi_m} I\{\xi_m \geq 1\}(X_{m,j} - \gamma_N) - \omega_m p(EY_m - \gamma_N), \quad m = 1,\ldots,N.$$

$$\gamma_N = \frac{1}{\Omega_N}\sum_{m=1}^{N}\omega_m EY_m, \quad \sigma_N^2 = p\sum_{m=1}^{N}\omega_m(E(Y_{mN} - \gamma_N)^2 - p(E(Y_{mN} - \gamma_N))^2). \tag{4.6}$$

From Theorems 4.3 and 4.4, we immediately have the following

**Corollary 4.3.** If (4.3) is satisfied, then for arbitrary $\delta \in (0,1]$ there exists a constant $C$ such that

$$\sup_{-\infty \leq u \leq \infty} |P\{S_{n,N} < u\sigma_N + n\gamma_N\} - \Phi(u)| \leq C\left(\beta_{2+\delta,N} + \frac{1}{(nq)^{\delta/2}}\right),$$

where

$$\beta_{2+\delta,N} = \sigma_N^{-(2+\delta)}\sum_{m=1}^{N} E|g_m(\xi_m)|^{2+\delta} \leq \frac{2^{2+\delta} p(1+p^{1+\delta})}{\sigma_N^{2+\delta}}\sum_{m=1}^{N}\omega_m^{2+\delta} E|Y_{m,N} - \gamma_N|^{2+\delta}. \tag{4.7}$$

**Corollary 4.4**. If (4.3) is satisfied, then there exist positive constants $c$ and $C$ such that

$$\Delta_N^{(5)} \leq C\left(\beta_{5,N} + \left(\frac{\bar{\omega}_N}{nq}\right)^{\frac{3}{2}} + \chi_N(c\tilde{T}_N, \beta_{5,N}^{-1})\right),$$

where $\tilde{T}_N = \min\left\{\beta_{3N}^{-1}, \sqrt{nq}/\bar{\omega}^2\right\}$, terms of the $\mathbb{W}_N^{(5)}(u,0)$ in (4.1) have the following forms:

$$\alpha_{1,2,N} = 0, \quad \alpha_{2,1,N} = \sqrt{\frac{q}{n}} \frac{\sum_{m=1}^{N}\omega_m(\alpha_{2,m} - 2p\alpha_{1,m}^2)}{\sum_{m=1}^{N}\omega_m(\alpha_{2,m} - p\alpha_{1,m}^2)}, \quad \alpha_{0,3,N} = \frac{1-2q}{\sqrt{nq}},$$

$$\alpha_{2,2,N} = \frac{\sum_{m=1}^{N}\omega_m(\alpha_{2,m}(1+(\omega_m-2)p) - \alpha_{1,m}^2(\omega_m-2)p(1-3q))}{\Omega_N p\sum_{m=1}^{N}\omega_m(\alpha_{2,m} - p\alpha_{1,m}^2)},$$

$$\sum_{m=1}^{N} \alpha_{11m}^2 = \frac{q\sum_{m=1}^{N} \omega_m^2 \alpha_{1,m}^2}{\Omega_N \sum_{m=1}^{N} \omega_m (\alpha_{2,m} - p\alpha_{1,m}^2)}, \qquad \sum_{m=1}^{N} \alpha_{20m}\alpha_{02m} = \frac{\sum_{m=1}^{N} \omega_m^2 (\alpha_{2,m} - p\alpha_{1,m}^2)}{\Omega_N \sum_{m=1}^{N} \omega_m (\alpha_{2,m} - p\alpha_{1,m}^2)},$$

$$\alpha_{3,0,N} = \sum_{m=1}^{N} \omega_m (\alpha_{3,m} - 3p\alpha_{1,m}\alpha_{2,m} - 2p^2 \alpha_{1,m}^3) \left( p^{2/3} \sum_{m=1}^{N} \omega_m (\alpha_{2,m} - p\alpha_{1,m}^2) \right)^{-3/2},$$

$$\alpha_{4,0,N} = \sum_{m=1}^{N} \omega_m (\alpha_{4,m} - 4p\alpha_{1,m}\alpha_{3,m} + 3(\omega_m - 1)p\alpha_{2,m}^2 - 6(\omega_m - 2)p^2\alpha_{1,m}^2 \alpha_{2,m}$$

$$- 3(3\omega_m - 2)p^3 \alpha_{1,m}^4) \left( \sqrt{p} \sum_{m=1}^{N} \omega_m (\alpha_{2,m} - p\alpha_{1,m}^2) \right)^{-2},$$

where $\alpha_{i,m} = E(Y_m - \gamma_N)^i$; also

$$\chi_N(c\tilde{T}_N, \beta_{5,N}^{-1}) \leq C\sqrt{n} \ln(\sigma_N^{-1} \beta_{5,N}^{-1}) \exp\left\{ -n\left(1 - \sup_{c\sigma_N \tilde{T}_N^{-1} \leq |t| \leq \sigma_N \beta_{5,N}^{-1}} \frac{1}{N} \sum_{m=1}^{N} |Ee^{itY_m}| \right) \right\}, \qquad (4.8)$$

and

$$\chi_N(c\tilde{T}_N, \beta_{5,N}^{-1}) \leq C\sqrt{n} \ln(\sigma_N^{-1} \beta_{5,N}^{-1}) \exp\left\{ -2nq\left(1 - \sup_{c\sigma_N \tilde{T}_N^{-1} \leq |t| \leq \sigma_N \beta_{5,N}^{-1}} \frac{1}{\Omega_N} \left| \sum_{m=1}^{N} \omega_m Ee^{itY_m} \right| \right) \right\}. \qquad (4.9)$$

In the case $\omega_1 = ... = \omega_N = 1$, our Corollary 4.3 improves a result of Mirakhmedov (1985), and a recent result of Zhao et al (2004) for the case when $(nq)^{-1/2} \leq \Delta_1$ (in their notations). Note that

$$\beta_{3N} = \frac{p}{\sigma_N^3} \sum_{m=1}^{N} E|Y_m - pEY_m - q\gamma_N|^3 + \frac{qp^3}{\sigma_N^3} \sum_{m=1}^{N} |E\,Y_m - \gamma_N|^3,$$

which provides a natural expression, showing the exact dependence of the bound on $p = n/N$, and moments of the elements of population, instead of the formula for $\Delta_2^*$ in Zhao et al (2004). This fact is confirmed by the second term in $\mathbb{W}_N^{(5)}(u,0)$ (see (4.1)), and that $\alpha_{1,2,N} = 0$. Also, in this case the three-term asymptotic expansion, i.e. $\mathbb{W}_N^{(5)}(u,0)$ coincides with that given by Mirakhmedov (1983); further, from our Corollary 4.4 follows the main result of Bloznelis (2000), and it extends Theorem 1 of Hu et al. (2007) giving an additional term in their asymptotic expansion for the case when $p$ is bounded away from one; this case is the most interesting in a sample scheme without replacement.

**4.3. Example C.** For this case we assume that $\mathcal{L}(\xi_m) = NB(d_m, p)$, with $p = n/(n + D_{1N})$, $m = 1,...,N$, where $D_{jN} = d_1^j + ... + d_N^j$, then $x_N = 0$. Putting $\rho = p/(1-p) = n/D_{1N}$ we get

$$B_N^2 = D_{1N}\rho(1+\rho), \quad \kappa_{4N} = (1 + 3\rho(1+\rho)(2 + D_{2N}D_{1N}^{-1}))(D_{1N}\rho(1+\rho))^{-1}.$$

**Theorem 4.5**. Let $D_{2N}D_{1N}^{-2} = o(N^{-1/2})$. If the Lindberg's condition (3.1) is satisfied then the DS $R_N(\eta)$ has asymptotic normal distribution with expectation $\Lambda_N$ and variance $\sigma_N^2$ as given in (2.1).

**Theorem 4.6**. There exists a constant $C$ such that $\Delta_N^{(3)} \leq C(\beta_{3N} + E_N)$, where

$$E_N = \frac{1}{\sqrt{n(1+\rho)}} + \sqrt{\frac{3}{D_{1N}}\left(2 + \frac{D_{2N}}{D_{1N}}\right)}\left(1 + \sqrt{\frac{3N}{D_{1N}}\left(2 + \frac{D_{2N}}{D_{1N}}\right)}\right).$$

**Remark 4.3**. Using the fact that $D_{1N}^2 \leq ND_{2N}$ we have

$$E_N \leq \frac{1}{\sqrt{n(1+\rho)}} + \frac{\sqrt{3}}{\sqrt{N}}\sqrt{1 + \frac{2N}{D_{1N}}}\left(1 + \sqrt{1 + \frac{2N}{D_{1N}}}\right).$$

Theorems 4.5 and 4.6 do improve as well as correct Theorems 13 and 14 of Mirakhmedov (1996).

**Theorem 4.7**. Let the statistic $R_N(\eta)$ be a lattice r.v. with span $h$ and a set of possible values $\Re$. If $E|g(\xi_m)|^5 < \infty$, $m = 1, 2, ..., N$ then uniformly in $z \in \Re$

$$\left|\frac{\sigma_N}{h} P\{R_N(\eta) = z\} - \frac{d}{du_z} W_N^{(5)}(u_z, 0)\right|$$

$$\leq C\left(\beta_{5,N} + \left(\frac{D_{3N}}{D_{1N}^3} + \frac{D_{2N}}{D_{1N}^2 n(1+\rho)}\right)^{3/4} + \sigma_N\sqrt{D_{1N}\rho(1+\rho)}\exp\left\{-\frac{1}{2}N\bar{H}_N\left(cT_N, \frac{\pi\sigma_N}{h}\right)\right\}\right),$$

where $u_z = (z - \Lambda_N)/\sigma_N$ .and $T_N$ is defined as in Sec 3 with $\mathcal{E}_N(1) = E_N$.

Now consider the following practical and important two-sample problem: Let $X_1, ..., X_{M-1}$ and $Y_1, ..., Y_n$ be two samples from continuous distributions $F$ and $G$ respectively defined on the same $A \subset R$. The classical two-sample problem is to test the null hypothesis of homogeneity $H_0: F = G$. Define the r.v.s

$$\eta_{m,k} = \sum_{i=1}^{n} I\{Y_i \in [X_{(m \cdot k)}, X_{(m \cdot k - k)}]\},$$

where $m = 1, ..., N$, $N = \lfloor M/k \rfloor$ is the largest integer that does not exceeds $M/k$, integer $k \geq 1$, $X_{(1)}, ..., X_{(M-1)}$ are the order statistics of the first sample $X_1, ..., X_{M-1}$. The r. vec. $(\eta_{1,k}, ..., \eta_{N,k})$ are called the "spacing-frequencies" i.e. frequencies of the second sample falling in between the spacings created by the first sample. A wide class of test statistics for testing $H_0$ can be expressed in the form (see Holst and Rao (1980) and Gatto and Jammalamadaka (1998))

$$V_N = \sum_{m=1}^{N} f_{m,N}(\eta_{m,k})$$

where $f_{m,N}$s are real- valued functions. It is easy to check that under $H_0$ the r.vec. $(\eta_{1,k}, ..., \eta_{N,k})$ satisfies equation (1.1) with $\mathcal{L}(\xi_m) = NB(k, p)$, $p = n/(n + M)$, i.e. the statistic $V_N$ is DS defined in the Pólya-Egenberger urn model. Hence, Theorems 4.5 and 4.6 immediately lead to the following Corollaries 4.5 and 4.6, by putting $d_m = k$, $m = 1,...,N$ and $\rho = n/M$.

**Corollary 4.5**. If the Lindeberg's condition (3.1) is satisfied then the statistic $V_N$ has asymptotic normal distribution with expectation $\Lambda_N$ and variance $\sigma_N^2$.

**Corollary 4.6.** There exists a constant $C > 0$ such that

$$\Delta_N^{(3)} \leq C\left(\beta_{3N} + \frac{1}{\sqrt{n(1+\rho)}} + \frac{1}{\sqrt{N}}\right).$$

Consider the following so-called Dixon's statistic: $\mathcal{D}_N = \sum_{m=1}^{N} \eta_{m,k}^2$. For this statistic we obtain:

$$\Lambda_N = M(1+(1+k)\rho), \quad \gamma_N = 1+2(1+k)\rho, \quad \sigma_N^2 = 2M(1+2k)\rho^2(1+\rho)^2,$$

$$g_m(\xi_m) = (\xi_m - k\rho)^2 - k\rho(1+\rho) - (1+2\rho)(\xi_m - k\rho),$$

$$\alpha_{12N} = \frac{\sqrt{2}(k+1)}{\sqrt{k(1+2k)}},$$

$$\alpha_{30N} = \frac{8k^2\rho^3(1+\rho)^3 + k(1+\rho)^2\,19 + 76\rho(1+\rho) + 2(1+\rho)^2(15-13\rho+16\rho^2(1+\rho))}{2\sqrt{2k}(1+2k)^{3/2}\rho^3(1+\rho)^3}.$$

Although the exact formula for $\beta_{4,N} = \alpha_{40N}$ is manageable, it is quite long and therefore we restrict ourselves to its leading term as $k \to \infty$, obtaining the following bounds

$$\beta_{4N} \leq CN^{-1}\max(1,\ k\rho(1+\rho)^{-4}), \qquad (\sigma_N \beta_{3N})^{-1} \geq c(k\rho(1+\rho))^{-1/2}.$$

The Dixon statistic satisfies the conditions of Theorem 4.7. In particular by evaluating the moments of the r.v. $g_m(\xi_m) = (\xi_m - k\rho)^2 - k\rho(1+\rho) - (1+2\rho)(\xi_m - k\rho)$, and using Corollaries 4.5, 4.6, and Theorem 4.7 along with the Remark 3.2, we obtain the following result. We omit the details.

**Corollary 4.7.**

(i) If $\sqrt{N}k^2\rho^2(1+\rho)^2 \to \infty$, then the Dixon statistic has an asymptotic normal distribution with mean $M(1+(1+k)\rho)$ and variance $2M(1+2k)\rho^2(1+\rho)^2$.

(ii) $P\{\mathcal{D}_N < u\sqrt{2M(1+2k)}\rho(1+\rho) + M(1+(1+k)\rho)\} = \Phi(u) + O\left(\frac{1}{\sqrt{N}} + \frac{1}{\sqrt{N}\,k^2\rho^2(1+\rho)^2}\right).$

(iii) Let $k \to \infty$ and $k = o(M^{1/3})$ then we have for any $b = 0, 1, ..., n(n-1)/2$

$$M(1+2k)\rho^2(1+\rho)^2 P(\mathcal{D}_N = n+2b) = \frac{e^{-u_b^2/2}}{\sqrt{2\pi}}\left(1 + \frac{u_b^2 - 3u_b}{6\sqrt{N}}\alpha_{30N} + \frac{(k+1)u}{\sqrt{2Nk(1+2k)}}\right) + O(N^{-1}),$$

where $u_b = (n+2b - M(1+(1+k)\rho))/\sqrt{2M(1+2k)}\rho(1+\rho)$.

Consider the class of symmetric tests (i.e. based on symmetric DS) for testing the hypothesis of homogeneity against some of "smooth" sequence of alternatives which approaches the null at the rate $O((nk)^{-1/4})$. The asymptotic power of symmetric tests increase as $k$ grows; the Dixon statistic is an example of a symmetric DS; it is known to be unique AMP within the class of symmetric tests for any fixed $k$, the step of spacings, see Jammalamadaka and Schweitzer (1985). Above stated Corollaries 4.4 -4.6 allow us to consider the situation when $k \to \infty$; in this case the AMP test is not unique. Comparison of the AMP tests based on their second order asymptotic efficiencies using the asymptotic expansion of the power function can be done. For this purpose, above presented asymptotic expansion results are central and such comparisons will be the subject of another investigation.

## 5. Appendix – Proofs:

**Proof of Proposition 2.1**. We need the following 3 Lemmas to complete the proof of this proposition.

**Lemma A.1**. Set $\ell_{s,N} = \min(\beta_{s,N}^{-1/s}, \kappa_{s,N}^{-1/s})$. There exist constants $c > 0$ and $C > 0$ such that if $\max |t|, |\tau| \leq c\ell_{s,N}$ then for $k = 0$ and $1$

$$\left| \frac{\partial^k}{\partial t^k} \left( \Psi_N(t,\tau) - e^{-\frac{t^2+\tau^2}{2}} \left( 1 + \sum_{v=1}^{s-3} N^{-v/2} P_v(t,\tau) \right) \right) \right| \leq C(\beta_{s,N} + \kappa_{s,N})(1 + |t|^s + |\tau|^s)e^{-\frac{t^2+\tau^2}{4}}.$$

Lemma A.1 follows from Theorem 9.11 of BR because of (2.2), and the fact that the sum of the r.v.s $(\hat{g}_m, \hat{\xi}_m)$ has unit correlation matrix. □

**Lemma A.2**. For any integer $l$ satisfying $0 \leq l \leq 3v$, where $v = 0, 1, ..., s-2$,

$$\left| \frac{\partial^l}{\partial t^l} G_{v,N}(t, x_N) \right| \leq c(l,v)(1 + (|t| + |x_N|)^{3v-l})(\beta_{v+2,N} + \kappa_{v+2,N}).$$

**Proof**. Similar to that of Lemma 9.5 of BR, p.71; the only difference being that in Equation (9.12) of BR, p.72, we use the inequality $\rho_{j_i+2} / \rho_2^{(j_i+2)/2} \leq (\rho_{r+2} / \rho_2^{(r+2)/2})^{j_i/r}$ to obtain

$$\left| \frac{\partial^l}{\partial t^l} P_{v,N}(t,\tau) \right| \leq c(l,v) \ 1 + (|t| + |\tau|)^{3v-l} \ (\beta_{v+2,N} + \kappa_{v+2,N})$$

$$\leq c(l,v) \ 1 + (|t| + |\tau|)^{3v-l} \ (\beta_{s,N}^{v/(s-2)} + \kappa_{s,N}^{v/(s-2)}). \tag{5.1}$$

Lemma A.2 follows from this and (2.6). □

**Lemma A.3.** Let $\max(\beta_{3N}, \kappa_{3N}) \leq 0.01$. If $|t| \leq 0.3 \beta_{3N}^{-1}$ and $|\tau| \leq 0.3 \kappa_{3N}^{-1}$ then for $k = 0$ and $1$

$$\left| \frac{\partial^k}{\partial t^k} \Psi_N(t,\tau) \right| \leq \exp\left\{ -\frac{t^2 + \tau^2}{10} \right\}.$$

Lemma A.3 follows from Lemma A, Part (2) of Mirakhmedov (2005). □

Put $T_N(s) = \min(\beta_{s,N}^{-1/s}, \kappa_{s,N}^{-1/s}, \mathcal{E}_N^{-1}(1))$, where $s \geq 3$. Note that $T_N(s) \leq T_N$, since $\beta_{s,N}^{-1/s} \leq \beta_{3,N}^{-1}$. Let $|t| \leq c_1 T_N(s)$, where $c_1 > 0$ is to be chosen sufficiently small. From (2.4) and (2.6)

$$\nabla_N(t) := \left| \frac{\partial^k}{\partial t^k} \left( \Theta_N(t, x_N) - e^{-\frac{t^2+x_N^2}{2}} \left( 1 + \sum_{v=1}^{s-3} N^{-v/2} G_{v,N}(t, x_N) \right) \right) \right|$$

$$\leq \int_{|\tau| \leq c_1 \ell_{s,N}} \left| \frac{\partial^k}{\partial t^k} \left( \Psi_N(t,\tau) - e^{-\frac{t^2+\tau^2}{2}} \sum_{v=1}^{s-3} N^{-v/2} P_{v,N}(t,\tau) \right) \right| d\tau$$

$$+ \int_{c_1 \ell_{s,N} \leq |\tau|} \left| \frac{\partial^k}{\partial t^k} \left( e^{-\frac{t^2+\tau^2}{2}} \sum_{v=1}^{s-3} N^{-v/2} P_{v,N}(t,\tau) \right) \right| d\tau + \int_{c_1 \ell_{s,N} \leq |\tau| \leq 0.3 \kappa_{3,N}^{-1}} \left| \frac{\partial^k}{\partial t^k} \Psi_N(t,\tau) \right| d\tau$$

$$+ \int_{0.3 \kappa_{3,N}^{-1} \leq |\tau| \leq \pi B_N} \left| \frac{\partial^k}{\partial t^k} \Psi_N(t,\tau) \right| d\tau = \mathfrak{I}_1 + \mathfrak{I}_2 + \mathfrak{I}_3 + \mathfrak{I}_4. \quad (5.2)$$

Applying Lemma A.1, (5.1) and Lemma A.3 to $\mathfrak{I}_1, \mathfrak{I}_2$ and $\mathfrak{I}_3$ respectively, we obtain, after some algebraic manipulation,

$$\mathfrak{I}_l \leq C(\beta_{s,N} + \kappa_{s,N})(1+|t|^s) e^{-\frac{t^2}{6}}, \quad l = 1, 2, 3. \quad (5.3)$$

Set $\hat{\psi}_m(t,\tau) = E\exp\{it\hat{g}_m(\xi_m) + i\tau\hat{\xi}_m\}$ and recall that $\Psi_N(t,\tau) = \prod_{m=1}^{N} \hat{\psi}_m(t,\tau)$. We have

$$|\hat{\psi}_m(t,\tau)|^2 = |\hat{\psi}_m(0,\tau)|^2 + E\left[ (e^{it\hat{g}_m^*} - 1)(e^{it\hat{\xi}_m^*} - 1) \right] + E(e^{it\hat{g}_m^*} - 1)$$

$$\leq |\hat{\psi}_m(0,\tau)|^2 + |t||\tau| E|\hat{g}_m^* \hat{\xi}_m^*| + t^2 E\hat{g}_m^{*2}$$

and

$$|\hat{\psi}_m(t,\tau)|^2 = |\hat{\psi}_m(0,\tau)|^2 + Ee^{it\hat{\xi}_m^*}(e^{it\hat{g}_m^*} - 1) \leq |\hat{\psi}_m(0,\tau)|^2 + 2|t| E|\hat{g}_m|.$$

Using these inequalities, and the fact that $x < \exp\{(x^2-1)/2\}$, and

$$\left| \frac{\partial}{\partial t} \prod_{m=1}^{N} \hat{\psi}_m(t,\tau) \right| \leq \sum_{m=1}^{N} \left| \frac{\partial}{\partial t} \hat{\psi}_m(t,\tau) \right| \prod_{l \neq m} |\hat{\psi}_l(t,\tau)|, \quad \sum_{m=1}^{N} \left| \frac{\partial}{\partial t} \hat{\psi}_m(t,\tau) \right| \leq |t| + 2|\tau|,$$

we find for k= 0, 1

$$\left| \frac{\partial^k}{\partial t^k} \prod_{m=1}^{N} \hat{\psi}_m(t,\tau) \right| \leq \sqrt{e}(|t|+|\tau|)^k \exp\left\{ \min(|t||\tau| + 2t^2, 2|t|\beta_{1,N}) - \frac{1}{2} \sum_{m=1}^{N} (1 - |\hat{\psi}_m(0,\tau)|^2) \right\}. \quad (5.4)$$

Choosing $c_1$ to be sufficiently small, using (5.4) and that $\beta_{1,N} \leq \sqrt{N}$, we get for $|t| \leq c_1 T_N(s)$

$$\mathfrak{I}_4 \leq CB_N^{k+1} \exp\left\{ -\frac{1}{4} M_N(0.3(B_N \kappa_{3,N})^{-1}) + \min(B_N|t| + 2t^2, 2|t|\sqrt{N}) \right\}$$

$$\leq CB_N^{k+1} \exp\left\{ -\frac{1}{8} M_N(0.3(B_N \kappa_{3,N})^{-1}) - \frac{t^2}{8} \right\}. \quad (5.5)$$

From (5.2),(5.3) and (5.4) it follows that

$$\left|\frac{\partial^k}{\partial t^k}\left(\Theta_N(t,x_N) - e^{-\frac{t^2+x_N^2}{2}}\left(1+\sum_{v=1}^{s-3} N^{-v/2} G_{v,N}(t,x_N)\right)\right)\right| \leq Ce^{-\frac{t^2}{8}}\Upsilon_{s,N}. \tag{5.6}$$

In particular (5.6) implies

$$\left|\Theta_N(0,x_N) - e^{-\frac{x_N^2}{2}}\sum_{v=0}^{s-3} N^{-v/2} G_{v,N}(0,x_N)\right| \leq C\Upsilon_{s,N}. \tag{5.7}$$

Put

$$\mathcal{G}_N^{(s)}(t,x_N) = e^{-\frac{t^2+x_N^2}{2}}\sum_{v=0}^{s-3} N^{-v/2} G_{v,N}(t,x_N) \ .$$

Then, from (2.5), (5.6) and (5.7) we have

$$\frac{\partial^k}{\partial t^k}(\varphi_N(t,x_N) - W_N^{(s)}(t,x_N)) = \frac{\partial^k}{\partial t^k}\left(\frac{\mathcal{G}_N^{(s)}(t,x_N)}{\mathcal{G}_N^{(s)}(0,x_N)} - W_N^{(s)}(t,x_N)\right)$$

$$+ \frac{\theta_2 \Upsilon_{s,N}}{\mathcal{G}_N^{(s)}(0,x_N) + \theta_2 \Upsilon_{s,N}}\left(\theta_1 t^k e^{-\frac{t^2}{8}} + \frac{1}{\mathcal{G}_N^{(s)}(0,x_N)}\frac{\partial^k}{\partial t^k}\mathcal{G}_N^{(s)}(t,x_N)\right), \tag{5.8}$$

where, as usual $|\theta_i| \leq 1$. Note that the polynomials $Q_{j,N}(x)$ in (2.7) are actually the results of the expansion of $\exp\left(x_N^2/2\right)\mathcal{G}_N^{(s)}(0,x_N)^{-1}$ taking into account that $G_{0,N}(x) = 1$; also it is obvious that $\mathcal{G}_N^{(s)}(0,x_N) \geq c$ for some $c > 0$. Using these facts and (2.7), (2.8), Lemma A.2 in (5.8) after some algebra we complete the proof of Proposition 2.1 for $|t| \leq c_1 T_N(s)$.

Let now $c_1 T_N(s) \leq |t| \leq c_1 T_N$. Then using Lemma A.2 it is easy to see that

$$\nabla_N(t) \leq \left|\frac{\partial^k}{\partial t^k}\left(\Theta_N(t,x_N) - e^{-\frac{t^2+x_N^2}{2}}\right)\right| + Ce^{-\frac{t^2}{8}}\Upsilon_{s,N}.$$

Apply above outlined technique for the first term in the rhs of this inequality using (2.4), Lemma A.1 with s=3, Lemma A.3 and the fact that $|t| \geq c_1 T_N(s)$, to complete the proof of Proposition 2.1; the details are omitted. □

**Proof of Theorem 3.1**. In addition to the notations of Sec.3, define

$$L_{1,N}(\varepsilon) = \frac{1}{N^{3/2}}\sum_{m=1}^{N} E|\hat{g}_m|^3 \mathbb{I}\{|\hat{g}_m| \leq \varepsilon\} \ .$$

Since $|x_N| \leq c_0$, Theorem 1 of Mirakhmedov (1996) gives: for arbitrary $\varepsilon > 0$ there exists a constant $C > 0$ such that

$$\Delta_N^{(3)} \leq C\left(L_{1,N}(\varepsilon) + L_{2,N}(\varepsilon) + \mathcal{L}_{1,N}(\varepsilon) + \mathcal{L}_{2,N}(\varepsilon) + B_N^2 \mathcal{L}_{1,N}(\varepsilon)\exp\left\{-\frac{1}{8}M_N(\pi(4B_N \mathcal{L}_{1,N}(\varepsilon))^{-1})\right\}\right.$$

$$+\frac{\max\left(\sqrt{M_N(\pi(4B_N\mathcal{L}_{1,N}(\varepsilon))^{-1})},\min(B_N,\sqrt{N})\right)}{M_N(\pi(4B_N\mathcal{L}_{1,N}(\varepsilon))^{-1})}\right).$$

Since $L_{1,N}(\varepsilon)\leq\varepsilon$, $\mathcal{L}_{1,N}(\varepsilon)\leq\varepsilon$ and $\varepsilon>0$ is arbitrarily small, Theorem 3.1 follows. $\square$

**Proof of Theorem 3.2.** From Theorem 2 of Mirakhmedov (1994)

$$\Delta_N^{(3)}\leq C\left(\beta_{2+\delta,N}+\kappa_{2+\delta,N}+B_N^2\kappa_{2+\delta,N}^{1/\delta}\exp\left\{-\frac{1}{8}M_N(\pi\ 4B_N\kappa_{2+\delta,N}^{1/\delta}\ ^{-1})\right\}+\mathcal{E}_N(\delta)\right),$$

since $P_N(u)=P(g_1(\xi_1)+...+g_N(\xi_N)<u\sigma_N|\zeta_N=n)$. If $M_N(\pi(4B_N\kappa_{2+\delta}^{1/\delta})^{-1})\leq cB_N$ for some $c>0$ then Theorem 3.2 is true with $C=c$. If $M_N(\pi(4B_N\kappa_{2+\delta,N}^{1/\delta})^{-1})>cB_N$ then

$B_N^2\kappa_{2+\delta,N}^{1/\delta}\exp\ -M_N(\pi(4B_N\kappa_{2+\delta,N}^{1/\delta})^{-1})/8\ \leq c_1\kappa_{2+\delta,N}^{1/\delta}\leq c_1\kappa_{2+\delta,N}$, since $\delta\in(0,1]$, and Theorem 3.2 follows. $\square$

**Proof of Theorem 3.3**. By the well-known Esseen's smoothing inequality we obtain

$$\nabla_N^{(s)}\leq\frac{1}{\pi}\int_{|t|\leq\beta_{s,N}^{-1}}\left|\frac{\varphi_N(t,x_N)-e^{-\frac{t^2}{2}}W_N^{(s)}(t,x_N)}{t}\right|dt+\frac{24}{\sqrt{2\pi}}\beta_{s,N}$$

$$\leq\frac{1}{\pi}\int_{|t|\leq cT_N}\left|\frac{\varphi_N(t,x_N)-e^{-\frac{t^2}{2}}W_N^{(s)}(t,x_N)}{t}\right|dt+\int_{cT_N\leq|t|\leq\beta_{s,N}^{-1}}\left|\frac{e^{-\frac{t^2}{2}}W_N^{(s)}(t,x_N)}{t}\right|dt$$

$$+\int_{cT_N\leq|t|\leq\beta_{s,N}^{-1}}\left|\frac{\varphi_N(t,x_N)}{t}\right|dt+\frac{24}{\sqrt{2\pi}}\beta_{s,N}.$$

We have

$$\left|\varphi_N(t,x_N)-e^{-\frac{t^2}{2}}W_N^{(s)}(t,x_N)\right|\leq|t|\max_{|u|\leq|t|}\left|\frac{\partial}{\partial u}\left(\varphi_N(u,x_N)-e^{-\frac{u^2}{2}}W_N^{(s)}(u,x_N)\right)\right|.$$

On the other hand, from the definition of $W_N^{(s)}(t,x_N)$, Lemma A.2, and the inequality $\beta_{3,N}\leq\beta_{s,N}^{1/(s-2)}$ we observe that

$$\int_{cT_N\leq|t|\leq\beta_{s,N}^{-1}}\left|\frac{e^{-\frac{t^2}{2}}W_N^{(s)}(t,x_N)}{t}\right|dt\leq Ce^{-c_3T_N^2}\leq c\Upsilon_{s,N}.$$

Therefore

$$\nabla_N^{(s)}\leq\frac{1}{\pi}\left[\int_{1\leq|t|\leq cT_N}\left|\varphi_N(t,x_N)-e^{-\frac{t^2}{2}}W_N^{(s)}(t,x_N)\right|dt\right.$$

$$\left.+\max_{|t|\leq1}\left|\frac{\partial}{\partial t}\left(\varphi_N(t,x_N)-e^{-\frac{t^2}{2}}W_N^{(s)}(t,x_N)\right)\right|\right]+c\Upsilon_{s,N}+\chi_N(cT_N,\beta_{s,N}^{-1}).$$

Applying Proposition 2.1 here completes the proof of Theorem 3.3. □

**Proof of Theorem 3.4** follows by standard methods outlined, as for instance in Petrov (1995, p.204-207), which uses the inversion formula and Propositions 2.1; the details are omitted. □

**Proof of Theorem 4.1**. The central moments of order $k$ of the Poi($\lambda$) r.v. is a polynomial in $\lambda$ of order $\lfloor k/2 \rfloor$, hence for even $j$ $\kappa_{j,N} \leq c(P_{(j-2)/2\,N} + n^{-(j-2)/2})$; for odd $j$ one can use the well-known inequality $\kappa_{l,N} \leq \kappa_{s,N}^{(l-2)/(s-2)}$, $3 \leq l \leq s$. Similarly, from the inequality (53) of Mirakhmedov (1996) we have $M_N(0.3\ B_N\kappa_{3N}^{-1}) \geq 0.2n(1+nP_{2N})^{-1}$. Theorem 4.1 follows from these facts and the inequality (3.7).

□

**Proof of Corollary 4.1**. To get the set of equalities given before Corollary 4.1, write $g_m(\xi_m) = (\xi_m - \lambda_m)^2 + 2\lambda_m(\xi_m - \lambda_m) - \lambda_m + (2nP_{2N} + 1)(\xi_m - \lambda_m)$, next to find a higher order central moments of the $Poi(\lambda)$ r.v. $\xi$ using the following recurrence formula of Kenney and Keeping (1953):

$$E(\xi - \lambda)^{\nu+1} = \nu\lambda E(\xi - \lambda)^{\nu-1} + \lambda \frac{d}{d\lambda} E(\xi - \lambda)^\nu.$$

Considering a r.v. which equals $p_m^{l-1}$ with the probability $p_m$, $m = 1,...,N$ and using well-known inequalities between moments one can see $P_{lN}^l \leq P_{l+1,N}^{l-1}$, $l = 2,3,...$; with equality iff $p_m = N^{-1}$, $m = 1,...,N$. Write $p_m = N^{-1}(1+\varepsilon_m)$, with $\varepsilon_m = Np_m - 1$, and put $\Sigma_N^2 = N^{-1}(\varepsilon_1^2 + ... + \varepsilon_N^2)$. It is easy to observe that $\sigma_N^2 = 2n\lambda(1 + c(1+\lambda)\Sigma_N^2)$, also $\alpha_{40N} \leq cn^2(1+\lambda^4\Sigma_N^2)/\sigma_N^4$. Considering separately the cases when $\lambda \to 0$, $\lambda \to \infty$, and $\lambda$ is bounded away from zero and infinity, one can show that $\beta_{4N} = N^{-1}\alpha_{40N} \leq c((n\lambda)^{-1} + N^{-1})$, $(\beta_{4N}\sigma_N^2)^{-1} \geq c(1+\lambda^2)^{-1}$ and $T_N \geq c(\sqrt{\lambda}\max(1,\sqrt{\lambda}))^{-1}$; the details omitted; here $c > 0$ is a constant and is different in different places. It is evident that the condition (3.9) is fulfilled. Corollary 4.1 follows from Theorem 4.1 and Remark 3.2. □

**Proof of Theorems 4.2-4.4**. We recall that in this case $\xi_m$ is $Bi\ \omega_m, p$ r.v. with $p = n/\Omega_N$. To find the central moments of the r.v. $\xi_m$ we use the following formula: for integer $k \geq 2$

$$E(\xi_m - \omega_m p)^k = pq\left(\frac{d}{dp}E(\xi_m - \omega_m p)^{k-1} + (k-1)\omega_m E(\xi_m - \omega_m p)^{k-2}\right). \tag{5.9}$$

We have: $B_N^2 = nq$,

$$\kappa_{2+\delta,N} \leq \kappa_{4,N}^{\delta/2} \leq \left(\frac{1-6pq}{nq} + 3\frac{\Omega_{2,N}}{\Omega_N^2}\right)^{\delta/2} \leq \left(\frac{1-6pq+3\overline{\omega}_N pq}{nq}\right)^{\delta/2} \leq \left(\frac{7\overline{\omega}_N}{4nq}\right)^{\delta/2},$$

$$\kappa_{5,N} \leq \kappa_{6,N}^{3/4} < 3\sqrt{3}\left(\frac{\overline{\omega}}{nq}\right)^{3/2}, \quad B_N\kappa_{3,N} \leq \sqrt{1-6pq+3\overline{\omega}_N pq} \leq \sqrt{7\overline{\omega}_N}.$$

Now use the inequalities

$$|E\exp\{i\tau\xi_m\}|^2 \leq \exp\{-4\omega_m pq\sin^2\tau/2\},$$

$$\sin^2\frac{\tau}{2} \geq \frac{\tau^2}{\pi^2}, |\tau| \leq \pi, \quad 1-e^{-u} \geq \frac{1-e^{-c}}{c}u, 0 \leq u \leq c, \quad (5.10)$$

to get

$$M_N(\pi(4B_N\kappa_{2+\delta,N}^{1/\delta})^{-1}) \geq \frac{(1-e^{-1})nq}{4\bar{\omega}_N(1-6pq+3nq\Omega_{2,N}\Omega_N^{-2})} \geq \frac{(1-e^{-1})nq}{4\bar{\omega}_N^2},$$

since $\Omega_{2,N} = \omega_1^2 + ... + \omega_N^2$. Finally $\mathcal{L}_{1,N}(\varepsilon) \leq \varepsilon^{-1}\kappa_{3,N} \leq \varepsilon^{-1}\sqrt{\bar{\omega}_N/nq}$. Theorems 4.2 - 4.4 follow from Theorems 3.1 - 3.3 respectively and the relations given above. □

**Proof of Corollaries 4.3 and 4.4**. Use inequality $(a_1 + ... + a_n)^s \leq n^{s-1}(a_1^s + ... + a_n^s)$, $a_m \geq 0, s \geq 1$, to get (4.7). Applying (5.9) we obtain the formulas for $\alpha_{ijN}$. Recall $\psi_m(t,\tau) = E\exp\{itf_{mN}(\xi_m) + i\tau\xi_m\}$ and fact that $\xi_m$ is a sum of $\omega_m$ independent $Bi(1,p)$ r.v.s. From (4.5) we have

$$|\psi_m(t,\tau)| = \left|\sum_{k=0}^{\omega_m} P(\xi_m = k)e^{i\tau k}Ee^{itf_{mN}(k)}\right| \leq \sum_{k=0}^{\omega_m} P(\xi_m = k)|Ee^{itY_m}|^k$$

$$\leq P(\xi_m = 0) + |Ee^{itY_m}|\sum_{k=1}^{\omega_m} P(\xi_m = k) = P(\xi_m = 0) + |Ee^{itY_m}|(1-P(\xi_m = 0)).$$

Hence

$$d_N \geq \frac{1}{N}\sum_{m=1}^{N}(1-(1-p)^{\omega_m})\left(1-\sup_{(\sigma_N T_N)^{-1} \leq |t| \leq \sigma_N \beta_{s,N}^{-1}}|Ee^{itY_m}|\right) \geq p\left(1-\sup_{(\sigma_N T_N)^{-1} \leq |t| \leq \sigma_N \beta_{s,N}^{-1}}\frac{1}{N}\sum_{m=1}^{N}|Ee^{itY_m}|\right),$$

since $P(\xi_m = 0) = (1-p)^{\omega_m}$. Inequality (4.8) follows from this and (3.5). On the other hand

$$\psi_m(t,\tau) = Ee^{(i\tau + \ln Ee^{itY_m})\xi_m} = (Ee^{i\tau\varsigma}(Ee^{itY_m})^\varsigma)^{\omega_m} = (1+p(Ee^{i(\tau+tY_m)}-1))^{\omega_m},$$

with $\mathcal{L}(\varsigma) = Bi(1,p)$. Hence

$$\prod_{m=1}^{N}|\psi_m(t,\tau)|^2 = \prod_{m=1}^{N}|1+p(Ee^{i(\tau+tY_m)}-1)|^{2\omega_m} \leq \exp\left\{-2pq\sum_{m=1}^{N}\omega_m(1-E\cos(\tau+tY_m))\right\}$$

$$\leq \exp\left\{-2\Omega_N pq\left(1-\frac{1}{\Omega_N}\left|\sum_{m=1}^{N}\omega_m Ee^{i\tau+itY_m}\right|\right)\right\}.$$

Inequality (4.9) follows. □

**Proof of Theorem 4.5 and 4.6**. Recall in this case $\mathcal{L}(\xi_m) = NB(d_m, p)$, with $p = n/(n+D_{1N})$, $m = 1,...,N$, where $D_{jN} = d_1^j + ... + d_N^j$, and $\rho = p/(1-p)$. We use that $Ee^{i\tau\xi_m} = (1-p)^{d_m}(1-pe^{i\tau})^{-d_m}$ to find the moments of the r.v. $\xi_m$ and that $B_N^2 = D_{1N}\rho(1+\rho)$,

$$B_N^2\kappa_{4N} = 1 + 3\rho^2(1+\rho)^2(2+D_{2N}D_{1N}^{-1}), \quad |Ee^{i\tau\xi_m}|^2 = \left(1+4\rho(1+\rho)\sin^2\frac{\tau}{2}\right)^{-d_m}.$$

Therefore, using the inequalities (5.11) we get

$$M_N(0.3(B_N \kappa_{3N})^{-1})^{-1}) \geq \frac{3(1-e^{-1/3})D_{1N}\rho(1+\rho)}{1+3\rho(1+\rho)(2+D_{2N}D_{1N}^{-1})} = 3(1-e^{-1/3})\kappa_{4N}^{-1}.$$

since $d_m \rho(1+\rho)(1+3\rho(1+\rho)(2+D_{2N}D_{1N}^{-1}))^{-1} < 1/3$. Therefore, Theorems 4.5, 4.6 and 4.7 follow from Theorems 3.1, 3.2 and 3.4 respectively and the inequality (3.7) by putting $\delta = 1$, and some simple algebra. □

---

[*] Mirakhmedov Sh. M. was formerly Mirakhmedov Sh. A